\documentstyle[amssymb,12pt]{article}

\newtheorem{th}{Theorem}[section]

\newtheorem{prop}[th]{Proposition}
\newtheorem{cor}[th]{Corollary}
\newtheorem{defn}[th]{Definition}
\newenvironment{defn-new}{\begin{defn} \em}{\end{defn}}
\newtheorem{rem}[th]{Remark}
\newenvironment{rem-new}{\begin{rem} \em}{\end{rem}}
\newtheorem{ex}[th]{Example}
\newenvironment{ex-new}{\begin{ex} \em}{\end{ex}}
\newtheorem{exer}[th]{Exercise}
\newenvironment{exer-new}{\begin{exer} \em}{\end{exer}}
\newtheorem{agr}[th]{Agreement}
\newenvironment{agr-new}{\begin{agr} \em}{\end{agr}}
\newtheorem{pbm}[th]{Problem}
\newenvironment{pbm-new}{\begin{pbm} \em}{\end{pbm}}

\makeatletter \@addtoreset{equation}{section} \makeatother

\begin{document}


\begin{center}
{\large {\bf Lightlike hypersurfaces of an \boldmath
$(\varepsilon)$-para Sasakian manifold}} \medskip

\bigskip

Selcen Y\"{u}ksel Perkta\c{s}$^{1}$, Erol K\i l\i \c{c}$^{2}$ and Mukut Mani
Tripathi$^{3}$\medskip

$^{1}$Department of Mathematics, Ad\i yaman University, Ad\i yaman, TURKEY

E-mail: sperktas@adiyaman.edu.tr

$^{2}$Department of Mathematics, \.{I}n\"{o}n\"{u} University, Malatya,
TURKEY

E-mail: erol.kilic@inonu.edu.tr

$^{3}$Department of Mathematics, Banaras Hindu University, Varanasi, INDIA

E-mail: mmtripathi66@yahoo.com

\medskip \medskip
\end{center}

\noindent {\bf Abstract. }In this paper, we initiate the study of lightlike
hypersurfaces of an $(\varepsilon )$-almost paracontact metric manifold
which are tangent to the structure vector field{\bf .} In particular, we
give definitions of invariant lightlike hypersurfaces and screen
semi-invariant lightlike hypersurfaces, and give some examples.
Integrability conditions for the distributions involved in the screen
semi-invariant lightlike hypersurface are investigated when the ambient
manifold is an $({\bf \varepsilon )}$-para Sasakian manifold.

\medskip


\medskip


\section{Introduction}

The theory of submanifolds of semi-Riemannian manifolds is one of the most
important topics in differential geometry. In case the induced metric on the
submanifold of semi-Riemannian manifold is degenerate, the study becomes
more difficult and is quite different from the study of nondegenerate
submanifolds. The primary difference between the lightlike submanifolds and
non-degenerate submanifolds arises due to the fact that in the first case
the normal vector bundle has non-trivial intersection with the tangent
vector bundle, and moreover in a lightlike hypersurface the normal vector
bundle is contained in the tangent vector bundle. Lightlike submanifolds of
semi-Riemannian manifolds were introduced by K.L.~Duggal and A.~Bejancu in
\cite{Dug-Bej-96}.

\medskip

In $1976$, an almost paracontact structure $(\phi ,\xi ,\eta )$ satisfying $%
\phi ^{2}=I-\eta \otimes \xi $ and $\eta (\xi )=1$ on a differentiable
manifold, was introduced by I.~S\={a}to \cite{Sato-76}. The structure is an
analogue of the almost contact structure \cite{Sasaki-60-Tohoku}, \cite%
{Blair-02-book} and is closely related to almost product structure (in
contrast to almost contact structure, which is related to almost complex
structure). An almost contact manifold is always odd-dimensional but an
almost paracontact manifold could be even-dimensional as well. In $1969$,
T.~Takahashi \cite{Takahashi-69-Tohoku-1} introduced almost contact
manifolds equipped with an associated pseudo-Riemannian metric. In
particular, he studied Sasakian manifolds equipped with an associated
pseudo-Riemannian metric. These indefinite almost contact metric manifolds
and indefinite Sasakian manifolds are also known as $\left( \varepsilon
\right) $-almost contact metric manifolds and $\left( \varepsilon \right) $%
-Sasakian manifolds, respectively (see \cite%
{Bej-Dug-93,Dug-90-IJMMS,Dug-Sah-07}). Lightlike hypersurfaces and
submanifolds of indefinite Sasakian manifolds were studied in $2003$ and $%
2007$ (see \cite{Kang-JKP-03} and \cite{Dug-Sah-07}). Also, in 1989,
K.~Matsumoto \cite{Mat-89} replaced the structure vector field $\xi $ by $%
-\, \xi $ in an almost paracontact manifold and associated a Lorentzian
metric with the resulting structure and called it a Lorentzian almost
paracontact manifold.

\medskip

An $(\varepsilon )$-Sasakian manifold is always odd-dimensional. On the
other hand, in a Lorentzian almost paracontact manifold given by
K.~Matsumoto, the semi-Riemannian metric has only index $1$ and the
structure vector field $\xi $ is always timelike. These circumstances
motivated the authors of \cite{Tri-KYK-10} to associate a semi-Riemannian
metric, not necessarily Lorentzian, with an almost paracontact structure,
and this indefinite almost paracontact metric structure is called an $\left(
\varepsilon \right) $-almost paracontact structure, where the structure
vector field $\xi $ is spacelike or timelike according as $\varepsilon =1$
or $\varepsilon =-1$ (see also \cite{Yuk-KTK-12}).

\medskip

In the present paper, as a first step to study lightlike geometry of $\left(
\varepsilon \right) $-almost paracontact metric manifolds we study lightlike
hypersurfaces. The paper is organized as follows. In section $2$, we give a
brief account of lightlike hypersurfaces of a semi-Riemannian manifold, for
later use. Section $3$\ is devoted to $\left( \varepsilon \right) $-almost
paracontact metric manifolds. In section $4$, we give investigate lightlike
hypersurfaces of an $\left( \varepsilon \right) $-almost paracontact metric
manifold. In section $5$, we define invariant lightlike hypersurfaces and
give an example. Screen semi invariant hypersurfaces are introduced in
Section $6$. Moreover, integrability conditions for the distributions
involved in the screen semi-invariant lightlike hypersurface are
investigated when the ambient manifold is an $(\varepsilon )$-para Sasakian
manifold.

\section{Lightlike Hypersurfaces\label{sect-lightlike}}

Let $(\widetilde{M},\widetilde{g})$ be an $(n+2)$-dimensional
semi-Riemannian manifold of fixed index $q\in \left \{ 1,\ldots
,n+1\right
\} $ and $M$ a hypersurface of $\widetilde{M}$. Assume that the
induced metric $g=\widetilde{g}|_{M}$ on the hypersurface is degenerate on $%
M.$ Then there exist a vector field $E\neq 0$ on $M$ such that
\[
g\left( E,X\right) =0,\qquad X\in \Gamma (TM).
\]%
The radical space \cite{ONeill-83} of $T_{x}M$, at each point $x\in M$, is
defined by
\begin{equation}
Rad~T_{x}M=\left \{ E\in T_{x}M:g\left( E,X\right) =0,\qquad X\in \Gamma
(T_{x}M)\right \} ,  \label{eq-null-2}
\end{equation}%
whose dimension is called the nullity degree of $g$ and $(M,g)$ is called a
lightlike hypersurface of $(\widetilde{M},\widetilde{g})$. Since $g$ is
degenerate and any null vector is perpendicular to itself, $T_{x}M^{\perp }$
is also degenerate and
\begin{equation}
Rad~T_{x}M=T_{x}M\cap T_{x}M^{\perp }.  \label{eq-null-3}
\end{equation}%
For a hypersurface $M$, $\dim T_{x}M^{\perp }=1$ implies that
\[
\dim Rad~T_{x}M=1,\qquad Rad~T_{x}M=T_{x}M^{\perp }.
\]%
We call $Rad~T\!M$ the radical distribution and it is spanned by the null
vector field $E$.

\medskip

Consider a complementary vector bundle $S(T\!M)$ of $Rad~T\!M$ in $T\!M$.
This means that
\begin{equation}
T\!M=S(T\!M)\perp Rad~T\!M,  \label{eq-null-4}
\end{equation}%
where $\perp $ denotes the orthogonal direct sum. The bundle $S(T\!M)$ is
called the screen distribution on $M$. Since the screen distribution $%
S(T\!M) $ is non-degenerate, there exists a complementary orthogonal vector
subbundle $S(T\!M)^{\perp }$ to $S(T\!M)$ in $T\widetilde{M}$ which is
called the screen transversal bundle. The rank of the screen transversal
bundle $S(T\!M)^{\perp }$ is $2$.

\medskip

Since $Rad~T\!M$ is a lightlike vector subbundle of $S(TM)^{\perp }$,
therefore for any local section $E\in \Gamma \left( Rad~T\!M\right) $ there
exists a unique local section $N$ of $S(T\!M)^{\perp }$ such that
\begin{equation}
\widetilde{g}(N,N)=0,\quad \quad \widetilde{g}(E,N)=1.  \label{eq-null-5}
\end{equation}%
Hence, $N$ is not tangent to $M$ and $\left \{ E,N\right \} $ is a local
frame field of $S(T\!M)^{\perp }$. Moreover, we have a $1$-dimensional
vector subbundle $ltr~T\!M$ of $T\widetilde{M}$, namely lightlike
transversal bundle, which is locally spanned by $N$. Then we set
\[
S(T\!M)^{\perp }=Rad~T\!M\oplus ltr~T\!M,
\]%
where the decomposition is not orthogonal. Thus we have the following
decomposition of $T\widetilde{M}$:
\begin{equation}
T\widetilde{M}=S(T\!M)\bot \,Rad~T\!M\oplus ltr~T\!M=T\!M\oplus ltr~T\!M.
\label{eq-null-6a}
\end{equation}%
From the above decomposition of a semi-Riemannian manifold $\widetilde{M}$
along a lightlike hypersurface $M$, we may consider the following local
quasi-orthonormal field of frames of $\widetilde{M}$ along $M$:
\[
\{X_{1},\ldots ,X_{n},E,N\},
\]%
where $\{X_{1},\ldots ,X_{n}\}$ is an orthonormal basis of $\Gamma (S(T\!M))$%
. According to the decomposition given by (\ref{eq-null-6a}), we have the
following Gauss and Weingarten formulas, respectively:
\begin{equation}
\widetilde{\nabla }_{X}Y=\nabla _{X}Y+B\left( X,Y\right) N,
\label{eq-null-7}
\end{equation}%
\begin{equation}
\widetilde{\nabla }_{X}N=-\,A_{N}X+\tau \left( X\right) N,  \label{eq-null-8}
\end{equation}%
where $B$ is a symmetric $\left( 0,2\right) $ tensor which is called the
second fundamental form and $A$ is an endomorphism of $T\!M$ which is called
the shape operator with respect to $N$ and $\tau $ is a $1$-form on $M$ \cite%
{Dug-Bej-96}. For each $X\in \Gamma (T\!M)$, we may write
\begin{equation}
X=PX+\theta \left( X\right) E,  \label{eq-null-9}
\end{equation}%
where $P$ is the projection of $T\!M$ on $S(T\!M)$ and $\theta $ is a $1$%
-form given by%
\begin{equation}
\theta \left( X\right) =\widetilde{g}\left( X,N\right) .  \label{eq-null-10}
\end{equation}%
From (\ref{eq-null-8}), for all $X,Y,Z\in \Gamma (T\!M)$, we get
\begin{equation}
\left( \nabla _{X}g\right) \left( Y,Z\right) =B(X,Y)\theta \left( Z\right)
+B\left( X,Z\right) \theta \left( Y\right) ,  \label{eq-null-11}
\end{equation}%
which implies that the induced connection $\nabla $ is a non-metric
connection on $M$. From (\ref{eq-null-4}), we have
\begin{equation}
\nabla _{X}W=\nabla _{X}^{\ast }W+C\left( X,W\right) E,  \label{eq-null-12}
\end{equation}%
\begin{equation}
\nabla _{X}E=-\,A_{E}^{\ast }X-\tau \left( X\right) E  \label{eq-null-13}
\end{equation}%
for all $X\in \Gamma (T\!M)$, $W\in \Gamma \left( S(T\!M)\right) $, where $C$%
, $A_{E}^{\ast }$ and $\nabla ^{\ast }$ are the local second fundamental
form, the local shape operator and the induced connection on $S(T\!M)$,
respectively. Note that $\nabla _{X}^{\ast }W$ and $A_{E}^{\ast }X$ belong
to $\Gamma \left( S(T\!M)\right) $. Also, we have the following identities
\begin{equation}
g\left( A_{E}^{\ast }X,W\right) =B\left( X,W\right) ,\quad g\left(
A_{E}^{\ast }X,N\right) =0,\quad B\left( X,E\right) =0,\quad g\left(
A_{N}X,N\right) =0.  \label{eq-null-14}
\end{equation}%
Moreover, from the first and third equations of (\ref{eq-null-14}) we have
\begin{equation}
A_{E}^{\ast }E=0.  \label{eq-null-15}
\end{equation}%
For more details we refer to \cite{Dug-Bej-96}, \cite{Dug-Jin-07} and \cite%
{Dug-Sah-10-book}.

\section{$(\protect\varepsilon )$-almost paracontact metric manifolds\label%
{sect-eps-ACMM}}

Let $\widetilde{M}$ be an almost paracontact manifold \cite{Sato-76}
equipped with an almost paracontact structure $(\phi ,\xi ,\eta )$
consisting of a tensor field $\phi $ of type $(1,1)$, a vector field $\xi $
and a $1$-form $\eta $ satisfying
\begin{equation}
\phi ^{2}=I-\eta \otimes \xi ,  \label{eq-phi-eta-xi}
\end{equation}%
\begin{equation}
\eta (\xi )=1,  \label{eq-eta-xi}
\end{equation}%
\begin{equation}
\phi \xi =0,  \label{eq-phi-xi}
\end{equation}%
\begin{equation}
\eta \circ \phi =0.  \label{eq-eta-phi}
\end{equation}%
Let $\widetilde{g}$ be a semi-Riemannian metric \cite{ONeill-83} such that
\begin{equation}
\widetilde{g}\left( \phi X,\phi Y\right) =\widetilde{g}\left( X,Y\right)
-\varepsilon \eta (X)\eta \left( Y\right) ,\qquad X,Y\in \Gamma (T\!%
\widetilde{M}),  \label{eq-metric-1}
\end{equation}%
where $\varepsilon =\pm 1$. Then $\widetilde{M}$ is called an $\left(
\varepsilon \right) $-almost paracontact metric manifold equipped with an $%
\left( \varepsilon \right) ${\em -}almost paracontact metric structure $%
(\phi ,\xi ,\eta ,\widetilde{g},\varepsilon )$ \cite{Tri-KYK-10}. In
particular, if ${\rm index}(g)=1$, that is when the metric is a Lorentzian
metric \cite{Beem-Ehrl-81}, then an $(\varepsilon )$-almost paracontact
metric manifold is called a Lorentzian almost paracontact manifold. From (%
\ref{eq-metric-1}) we have
\begin{equation}
\widetilde{g}\left( X,\phi Y\right) =\widetilde{g}\left( \phi X,Y\right)
\label{eq-metric-2}
\end{equation}%
along with
\begin{equation}
\widetilde{g}\left( X,\xi \right) =\varepsilon \eta (X)  \label{eq-metric-3}
\end{equation}%
for all $X,Y\in \Gamma (T\! \widetilde{M})$. From (\ref{eq-metric-3}) it
follows that
\begin{equation}
\widetilde{g}\left( \xi ,\xi \right) =\varepsilon ,  \label{eq-g(xi,xi)}
\end{equation}%
that is, the structure vector field $\xi $ is never lightlike.

\medskip

Let $(\widetilde{M},\phi ,\xi ,\eta ,\widetilde{g},\varepsilon )$ be an $%
(\varepsilon )$-almost paracontact metric manifold (resp. a Lorentzian
almost paracontact manifold). If $\varepsilon =1$, then $\widetilde{M}$ will
be said to be a spacelike $(\varepsilon )$-almost paracontact metric
manifold (resp. a spacelike Lorentzian almost paracontact manifold).
Similarly, if $\varepsilon =-\,1$, then $\widetilde{M}$ will be said to be a
timelike $(\varepsilon )$-almost paracontact metric manifold (resp. a
timelike Lorentzian almost paracontact manifold) \cite{Tri-KYK-10}.

\medskip

An $\left( \varepsilon \right) $-almost contact metric structure is called
an $\left( \varepsilon \right) ${\em -}para Sasakian structure if
\begin{equation}
(\widetilde{\nabla }_{X}\phi )Y=-\, \widetilde{g}(\phi X,\phi Y)\xi
-\varepsilon \eta \left( Y\right) \phi ^{2}X,\qquad X,Y\in \Gamma (T\!%
\widetilde{M}),  \label{para2}
\end{equation}%
where $\widetilde{\nabla }$ is the Levi-Civita connection with respect to $%
\widetilde{g}$. A manifold endowed with an $\left( \varepsilon \right) $%
-para Sasakian structure is called an $\left( \varepsilon \right) ${\em -}%
para Sasakian manifold \cite{Tri-KYK-10}. In an $\left( \varepsilon \right) $%
{\em -}para Sasakian manifold, we have
\begin{equation}
\widetilde{\nabla }\xi =\varepsilon \phi  \label{para3}
\end{equation}%
and
\begin{equation}
\Phi \left( X,Y\right) =\widetilde{g}\left( \phi X,Y\right) =\varepsilon
\widetilde{g}(\widetilde{\nabla }_{X}\xi ,Y)=(\widetilde{\nabla }_{X}\eta
)Y,~X,Y\in \Gamma (T\widetilde{M})  \label{para4}
\end{equation}%
where
\begin{equation}
\Phi \left( X,Y\right) =\widetilde{g}\left( X,\phi Y\right) .  \label{para5}
\end{equation}%
From (\ref{para5}) we have
\begin{equation}
\Phi \left( X,\xi \right) =0.  \label{para6}
\end{equation}

In an $\left( \varepsilon \right) $-para Sasakian manifold the following
equations hold for any $X,Y,Z\in \Gamma (T\widetilde{M})$ \cite{Tri-KYK-10}:

\begin{equation}
\widetilde{R}\left( X,Y\right) \xi =\eta \left( X\right) Y-\eta \left(
Y\right) X,  \label{para7}
\end{equation}%
\begin{equation}
\widetilde{R}\left( X,Y,Z,\xi \right) =-\, \eta \left( X\right) \widetilde{g}%
\left( Y,Z\right) +\eta \left( Y\right) \widetilde{g}\left( X,Z\right) ,
\label{para8}
\end{equation}%
\begin{equation}
\eta (\widetilde{R}\left( X,Y\right) Z)=-\, \varepsilon \eta \left( X\right)
\widetilde{g}\left( Y,Z\right) +\varepsilon \eta \left( Y\right) \widetilde{g%
}\left( X,Z\right) ,  \label{para9}
\end{equation}%
\begin{equation}
\widetilde{R}\left( \xi ,X\right) Y=-\, \varepsilon \widetilde{g}\left(
X,Y\right) \xi +\eta \left( Y\right) X,  \label{para10}
\end{equation}%
\begin{equation}
\widetilde{S}\left( Y,\xi \right) =-\left( n-1\right) \eta \left( Y\right) ,
\label{para11}
\end{equation}%
where $\widetilde{R}$ is the Riemannian curvature tensor and $S$ is the
Ricci tensor of $\widetilde{M}$.

\section{Lightlike hypersurfaces of $\left( \protect\varepsilon \right) $%
-para Sasakian manifolds\label{sect-main-results}}

Let $(\widetilde{M},\phi ,\xi ,\eta ,\widetilde{g},\varepsilon )$ be an $%
\left( n+2\right) $-dimensional $(\varepsilon )$-para Sasakian manifold and $%
M$ be a lightlike hypersurface of $\widetilde{M}$, such that the structure
vector field $\xi $ is tangent to $M$. Since $\xi $ is a non-null vector
field, it belongs to the screen distribution $S(T\!M)$. If index$\left(
\widetilde{g}\right) =1$, in order that $M$ is a lightlike hypersurface, it
is necessary that the structure vector field $\xi $ must be a spacelike
vector field, that is, $\widetilde{M}$ must be a spacelike para Sasakian
manifold. If {\rm index}$\left( \widetilde{g}\right) >1$, then $\widetilde{M}
$ may also be a timelike para Sasakian manifold.

\medskip

For local sections $E$ and $N$ of $Rad~T\!M$ and $ltr~T\!M$, respectively,
in view of (\ref{eq-metric-3}), we have
\[
\eta \left( E\right) =0,\qquad \eta \left( N\right) =0.
\]%
From (\ref{eq-metric-1}), it is easy to see that $\phi E$ and $\phi N$ are
lightlike vector fields and
\[
\phi ^{2}E=E,\qquad \phi ^{2}N=N.
\]%
Now, for $X\in \Gamma (T\!M)$, we write
\begin{equation}
\phi X=\varphi X+u\left( X\right) N,  \label{eq-llh-1}
\end{equation}%
where $\varphi X\in \Gamma (T\!M)$ and
\begin{equation}
u\left( X\right) =\widetilde{g}\left( \phi X,E\right) =\widetilde{g}\left(
X,\phi E\right) .  \label{eq-llh-2}
\end{equation}

\begin{prop}
Let $(\widetilde{M},\phi ,\xi ,\eta ,\widetilde{g},\varepsilon )$ be an $%
\left( n+2\right) $-dimensional $(\varepsilon )$-para Sasakian manifold and $%
M$ be a lightlike hypersurface of $\widetilde{M}$, such that the structure
vector field $\xi $ is tangent to $M$. Then we have
\begin{equation}
\widetilde{g}\left( \phi E,E\right) =0,  \label{eq-g(phiE,E)=0}
\end{equation}%
\begin{equation}
\widetilde{g}\left( \phi E,N\right) =\varepsilon g\left( A_{N}E,\xi \right) ,
\label{eq-g(phiE,N)}
\end{equation}%
where $E$ is a local section of $Rad~T\!M$ and $N$ is a local section of $%
ltr~T\!M$.
\end{prop}

\noindent {\bf Proof.} From (\ref{para3}) and (\ref{eq-null-15}), we get (%
\ref{eq-g(phiE,E)=0}). By using (\ref{para3}), (\ref{eq-null-5}) and (\ref%
{eq-null-8}) we obtain
\[
\varepsilon \widetilde{g}\left( \phi E,N\right) =\widetilde{g}(\widetilde{%
\nabla }_{E}\xi ,N)=-\widetilde{g}(\xi ,\widetilde{\nabla }_{E}N)=g\left(
A_{N}E,\xi \right) ,
\]%
which implies (\ref{eq-g(phiE,N)}). $\blacksquare $

\medskip

\begin{rem-new}
From (\ref{eq-g(phiE,E)=0}) we see that there is no component of $\phi E$ in
$ltr~T\!M$, thus $\phi E\in \Gamma (T\!M)$. Moreover, (\ref{eq-g(phiE,N)})
implies that there may be a component of $\phi E$ in $Rad~T\!M$. Thus, in
view of (\ref{eq-null-9}) and (\ref{eq-g(phiE,E)=0}), we observe that
\begin{equation}
\phi E=\varphi E=P\phi E+\theta \left( \phi E\right) E.  \label{Light-3'-}
\end{equation}
\end{rem-new}

\begin{prop}
Let $(\widetilde{M},\phi ,\xi ,\eta ,\widetilde{g},\varepsilon )$ be an $%
\left( n+2\right) $-dimensional $(\varepsilon )$-almost paracontact metric
manifold and $M$ be a lightlike hypersurface of $\widetilde{M}$ such that
the structure vector field $\xi $ is tangent to $M$. Then we have
\begin{equation}
g\left( X,\varphi Y\right) =g\left( \varphi X,Y\right) +(u\wedge \theta
)(X,Y),  \label{eq-g(X,phiY)}
\end{equation}%
\begin{eqnarray}
g\left( \varphi X,\varphi Y\right) &=&g\left( X,Y\right) -\varepsilon \eta
(X)\eta \left( Y\right)  \nonumber \\
&&-\,u(X)\theta (\varphi Y)-u(Y)\theta (\varphi X)  \label{eq-g(phiX,phiY)}
\end{eqnarray}%
for any $X,Y\in \Gamma (T\!M)$
\end{prop}

\noindent {\bf Proof.} From (\ref{eq-llh-1}) and (\ref{eq-llh-2}), we get
\[
\widetilde{g}\left( \phi X,Y\right) =g\left( \varphi X,Y\right) +u(X)\theta
(Y).
\]%
Hence in view of (\ref{eq-metric-2}) we get (\ref{eq-g(X,phiY)}). Using (\ref%
{eq-llh-1}) we have%
\begin{equation}
\widetilde{g}\left( \phi X,\phi Y\right) =g\left( \varphi X,\varphi Y\right)
+u(X)\theta (\varphi Y)+u(Y)\theta (\varphi X).  \label{eq-g(vphiX,vphiY)}
\end{equation}%
Thus by using (\ref{eq-g(vphiX,vphiY)}) and (\ref{eq-metric-1}) we complete
the proof. $\blacksquare $

\medskip

\begin{cor}
Let $(\widetilde{M},\phi ,\xi ,\eta ,\widetilde{g},\varepsilon )$ be an $%
\left( n+2\right) $-dimensional $(\varepsilon )$-almost paracontact metric
manifold and $M$ be a lightlike hypersurface of $\widetilde{M}$, such that
the structure vector field $\xi $ is tangent to $M$. Then we have
\[
g\left( \xi ,\varphi X\right) =0,\qquad X\in \Gamma (T\!M).
\]
\end{cor}

\begin{prop}
Let $(\widetilde{M},\phi ,\xi ,\eta ,\widetilde{g},\varepsilon )$ be an $%
\left( n+2\right) $-dimensional $(\varepsilon )$-para Sasakian manifold and $%
M$ be a lightlike hypersurface of $\widetilde{M}$, such that the structure
vector field $\xi $ is tangent to $M$. Then, for any $X\in \Gamma (T\!M)$ we
have
\begin{equation}
\varphi ^{2}X=X-\eta (X)\xi -u(\varphi X)N-u(X)\phi N,  \label{eq-phi-sq-X-1}
\end{equation}%
\begin{equation}
\varphi X=\varepsilon \nabla _{X}\xi ,  \label{eq-phi-X-1}
\end{equation}%
\begin{equation}
B(X,\xi )=\varepsilon u(X).  \label{eq-B(X,xi)}
\end{equation}
\end{prop}

\noindent {\bf Proof.} From (\ref{eq-llh-1}) and (\ref{eq-phi-eta-xi}), we
get (\ref{eq-phi-sq-X-1}). Next, from (\ref{para3}), (\ref{eq-null-7}) and (%
\ref{eq-llh-1}) we have
\[
\varepsilon \nabla _{X}\xi +\varepsilon B(X,\xi )N=\varphi X+u\left(
X\right) N.
\]%
Then by equating the tangential and the transversal parts in the previous
equation we get (\ref{eq-phi-X-1}) and (\ref{eq-B(X,xi)}), respectively. $%
\blacksquare $

\section{Invariant lightlike hypersurfaces\label{sect-inv-lightlike-hyp}}

We begin with the following

\begin{defn-new}
Let $(\widetilde{M},\phi ,\xi ,\eta ,\widetilde{g},\varepsilon )$ be a $%
\left( n+2\right) $-dimensional an $\left( \varepsilon \right) $-almost
paracontact metric manifold and $M$ be a lightlike hypersurface of $%
\widetilde{M}$. If $\phi \left( S(T\!M)\right) =S(T\!M)$, then $M$ will be
called an invariant lightlike hypersurface of $\widetilde{M}$.
\end{defn-new}

\begin{ex-new}
Let ${\Bbb R}^{5}$\ be the $5$-dimensional real number space with a
coordinate system $\left( x,y,z,t,s\right) $. Defining
\[
\eta =ds-ydx-tdz\ ,\qquad \xi =\frac{\partial }{\partial s}\ ,
\]%
\[
\phi \left( \frac{\partial }{\partial x}\right) =-\, \frac{\partial }{%
\partial x}-y\frac{\partial }{\partial s}\ ,\qquad \phi \left( \frac{%
\partial }{\partial y}\right) =-\, \frac{\partial }{\partial y}\ ,
\]%
\[
\phi \left( \frac{\partial }{\partial z}\right) =-\, \frac{\partial }{%
\partial z}-t\frac{\partial }{\partial s}\ ,\qquad \phi \left( \frac{%
\partial }{\partial t}\right) =-\, \frac{\partial }{\partial t}\ ,\qquad
\phi \left( \frac{\partial }{\partial s}\right) =0\ ,
\]%
\begin{eqnarray*}
\widetilde{g} &=&-\, \left( dx\right) ^{2}-\left( dy\right) ^{2}+\left(
dz\right) ^{2}+\left( dt\right) ^{2}+\left( ds\right) ^{2} \\
&&-\,t\left( dz\otimes ds+ds\otimes dz\right) -y\left( dx\otimes
ds+ds\otimes dx\right) ,
\end{eqnarray*}%
the set $(\phi ,\xi ,\eta ,g)$\ is a spacelike $\left( \varepsilon \right) $%
-almost paracontact structure with {\rm index}$\left( g\right) =3$ on ${\Bbb %
R}^{5}$. Consider a hypersurface $M$ of ${\Bbb R}^{5}$ given by $y=t$. It is
easy to check that $M$ is a lightlike hypersurface whose radical
distribution $Rad~T\!M$ is spanned by
\[
E=\frac{\partial }{\partial y}+\frac{\partial }{\partial t}.
\]%
Then the lightlike transversal vector bundle $ltr~T\!M$ is spanned by
\[
N=\frac{1}{2}\left( -\frac{\partial }{\partial y}+\frac{\partial }{\partial t%
}\right) ,
\]%
and the screen bundle $S(T\!M)$ is spanned by%
\[
\left \{ U_{1},U_{2},\xi \right \} ,
\]%
where $U_{1}=\frac{\partial }{\partial x}$ and $U_{2}=\frac{\partial }{%
\partial z}$. We easily check that
\[
\phi E=-E,\quad \phi N=-N.
\]%
Thus $M$ is a invariant lightlike hypersurface of ${\Bbb R}^{5}$.
\end{ex-new}

In the following we give a characterization of an invariant lightlike
hypersurface.

\begin{th}
Let $(\widetilde{M},\phi ,\xi ,\eta ,\widetilde{g},\varepsilon )$ be an $%
(\varepsilon )$-almost paracontact metric manifold. Then $M$ is an invariant
lightlike hypersurface of $\widetilde{M}$ if and only if
\[
\phi Rad~T\!M=Rad~T\!M\qquad {\rm and}\qquad \phi \,ltr~T\!M=ltr~T\!M.
\]
\end{th}

\noindent {\bf Proof.} Let $M$ be an invariant lightlike hypersurface of $%
\widetilde{M}$. From (\ref{Light-3'-}), for any $X\in \Gamma (T\!M)$, we get
$g\left( P\phi E,PX\right) =0$, that is, there is no component of $\phi E$
in $S(T\!M)$ and $\phi Rad~T\!M=Rad~T\!M$. For any local section $N$ of $%
ltr~T\!M$, we can write
\begin{equation}
\phi N=P\phi N+\widetilde{g}\left( \phi N,N\right) E+\widetilde{g}\left(
\phi N,E\right) N.  \label{Light-3''-}
\end{equation}%
From (\ref{Light-3''-}), for any $X\in \Gamma (T\!M)$, we get $g\left( P\phi
N,PX\right) =0$, that is, there is no component of $\phi N$ in $S(T\!M)$. If
we apply $\phi $ to (\ref{Light-3''-}), then we get
\[
2\widetilde{g}\left( \phi N,N\right) \widetilde{g}\left( \phi N,E\right) =0.
\]%
Since $\ker \phi ={\rm Span}\left \{ \xi \right \} $, we obtain $\widetilde{g%
}\left( \phi N,N\right) =0$. Thus we get $\phi N=\widetilde{g}\left( \phi
N,E\right) N$, that is $\phi \,ltr~T\!M=ltr~T\!M$.

\smallskip

Conversely, let $\phi Rad~T\!M=Rad~T\!M$ and $\phi \,ltr~T\!M=ltr~T\!M$. For
any $X\in \Gamma (S(T\!M))$ we have
\[
\widetilde{g}\left( \phi X,E\right) =\widetilde{g}\left( X,\phi E\right) =0;
\]%
thus there is no component of $\phi X$ in $ltr~T\!M$ . Similarly, we get
\[
\widetilde{g}\left( \phi X,N\right) =\widetilde{g}\left( X,\phi N\right) =0,
\]%
which implies that there is no component of $\phi X$ in $Rad~T\!M$. This
completes the proof. $\blacksquare $ \medskip

\begin{th}
\label{th-llh-inv-1} Let $(\widetilde{M},\phi ,\xi ,\eta ,\widetilde{g}%
,\varepsilon )$ be an $(\varepsilon )$-almost paracontact metric manifold.
Let $M$ be an invariant lightlike hypersurface of $\widetilde{M}$. Then $%
(M,\varphi ,\xi ,\eta ,g,\varepsilon )$ is an $(\varepsilon )$-almost
paracontact metric manifold.
\end{th}

\noindent {\bf Proof.} Let $M$ be an invariant lightlike hypersurface of $%
\widetilde{M}$. Let us assume that $X,Y\in \Gamma (T\!M)$. From (\ref%
{eq-llh-1}), we get
\begin{equation}
\phi X=\varphi X.  \label{eq-phiX=}
\end{equation}%
Using (\ref{eq-phi-eta-xi}) and (\ref{eq-phiX=}), we have
\begin{equation}
\varphi ^{2}X=X-\eta (X)\xi .  \label{eq-phi-sq-X-inv}
\end{equation}%
Also from (\ref{eq-phiX=}), it follows that
\begin{equation}
\varphi \xi =0.  \label{eq-phi-xi=0-inv}
\end{equation}%
Next, in view of (\ref{eq-phi-sq-X-inv}) and (\ref{eq-phi-xi=0-inv}) one can
easily see that
\[
\eta \circ \varphi =0,
\]%
\[
\eta \left( \xi \right) =1.
\]%
Moreover, from (\ref{eq-g(phiX,phiY)}) we have
\[
g\left( \varphi X,\varphi Y\right) =g\left( X,Y\right) -\varepsilon \eta
(X)\eta \left( Y\right) .
\]%
This completes the proof. $\blacksquare $

\begin{prop}
Let $M$ be an invariant lightlike hypersurface of an $(\varepsilon )$-para
Sasakian manifold $(\widetilde{M},\phi ,\xi ,\eta ,\widetilde{g},\varepsilon
)$. Then we have
\[
g\left( A_{N}PX,\xi \right) =0,\qquad X\in \Gamma \left( TM\right) .
\]
\end{prop}

\noindent {\bf Proof. }Since $g\left( \xi ,N\right) =0$, using (\ref{para3}%
), we get
\[
\widetilde{g}\left( \widetilde{\nabla }_{X}N,\xi \right) =-\varepsilon
\widetilde{g}\left( N,\phi X\right) .
\]%
From (\ref{eq-null-8}), we have the assertion of the proposition. $%
\blacksquare $

\begin{th}
\label{th-llh-inv-2} An invariant lightlike hypersurface of an $(\varepsilon
)$-para Sasakian manifold is always $(\varepsilon )$-para Sasakian.
Moreover,
\begin{equation}
B\left( X,\varphi Y\right) N-B\left( X,Y\right) \phi N=0,
\label{eq-llh-inv-B}
\end{equation}%
\begin{equation}
\varphi \left( A_{N}X\right) =A_{\phi N}X-\theta \left( X\right) \xi
\label{eq-phi-A-N-X}
\end{equation}%
for all $X,\,Y\in {\Gamma ({TM})}$.
\end{th}

\noindent {\bf Proof.} We have%
\begin{eqnarray*}
(\widetilde{\nabla }_{X}\phi )Y &=&\widetilde{\nabla }_{X}\phi Y-\phi (%
\widetilde{\nabla }_{X}Y) \\
&=&\widetilde{\nabla }_{X}\phi Y-\phi (\nabla _{X}Y+B\left( X,Y\right) N) \\
&=&\nabla _{X}\varphi Y+B\left( X,\varphi Y\right) N-\varphi \nabla
_{X}Y-B\left( X,Y\right) \phi N \\
&=&\left( \nabla _{X}\varphi \right) Y+B\left( X,\varphi Y\right) N-B\left(
X,Y\right) \phi N,
\end{eqnarray*}%
which in view of (\ref{para2}) gives
\begin{equation}
-\,g(\varphi X,\varphi Y)\xi -\varepsilon \eta \left( Y\right) \varphi
^{2}X=\left( \nabla _{X}\varphi \right) Y+B\left( X,\varphi Y\right)
N-B\left( X,Y\right) \phi N.  \label{eq-med-1}
\end{equation}%
Equating tangential parts in (\ref{eq-med-1}) provides
\begin{equation}
\left( \nabla _{X}\varphi \right) Y=-\,g(\varphi X,\varphi Y)\xi
-\varepsilon \eta \left( Y\right) \varphi ^{2}X,  \label{eq-llh-inv-PS}
\end{equation}%
In view of (\ref{eq-llh-inv-PS}) and Theorem~\ref{th-llh-inv-1} we see that $%
M$ is $(\varepsilon )$-para Sasakian. Equating transversal parts in (\ref%
{eq-med-1}) yields (\ref{eq-llh-inv-B}).

\smallskip

Next, using (\ref{para2}) and (\ref{eq-null-8}) we have
\begin{eqnarray*}
-\theta \left( X\right) \xi &=&(\widetilde{\nabla }_{X}\phi )N=\widetilde{%
\nabla }_{X}\phi N-\phi (\widetilde{\nabla }_{X}N) \\
&=&-A_{\phi N}X+\widetilde{g}\left( \widetilde{\nabla }_{X}\phi N,E\right)
N+\phi \left( A_{N}X\right) -\tau \left( X\right) \phi N.
\end{eqnarray*}%
In the last equation, if we equate the tangential parts, we get
\[
-\theta \left( X\right) \xi =-\,A_{\phi N}X+\varphi (A_{N}X).
\]%
This completes the proof. $\blacksquare $

\begin{rem}
It is well-known that, if there exists a lightlike hypersurface in an $%
(\varepsilon )$- Sasakian manifold, then \ the dimension of the Sasakian
manifold must be greater than or equal to $5$. But in the case of an $%
(\varepsilon )$-paracontact metric manifold there is no such restriction on
the dimension of the ambient manifold for the existence of lightlike
hypersurfaces.
\end{rem}

\section{Screen semi-invariant lightlike hypersurfaces\label{sect-SSILH}}

We begine with the following:

\begin{defn-new}
Let $(\widetilde{M},\phi ,\xi ,\eta ,\widetilde{g},\varepsilon )$ be a $%
\left( n+2\right) $-dimensional $\left( \varepsilon \right) $-almost
paracontact metric manifold and $M$ be a lightlike hypersurface of $%
\widetilde{M}$. If $\phi Rad~T\!M\subset S(T\!M)$ and $\phi
\,ltr~T\!M\subset S(T\!M)$, then $M$ will be called a screen semi-invariant
lightlike hypersurface of $\widetilde{M}$.
\end{defn-new}

\begin{ex-new}
Let ${\Bbb R}^{5}$\ be the $5$-dimensional real number space with a
coordinate system $\left( x,y,z,t,s\right) $. We define
\[
\eta =\frac{1}{2}\left( zdx+tdy+ds\right) ,\qquad \xi =2\frac{\partial }{%
\partial s},
\]%
\[
\phi X=-X_{3}\frac{\partial }{\partial x}-X_{4}\frac{\partial }{\partial y}%
-X_{1}\frac{\partial }{\partial z}-X_{2}\frac{\partial }{\partial t}+\left(
X_{3}z+X_{4}t\right) \frac{\partial }{\partial s},
\]%
\[
\widetilde{g}=-\frac{1}{4}\left( dx\otimes dx+dz\otimes dz\right) +\frac{1}{4%
}\left( dy\otimes dy+dt\otimes dt\right) +\eta \otimes \eta .
\]%
Here $X$ is a vector field given by%
\[
X=X_{1}\frac{\partial }{\partial x}+X_{2}\frac{\partial }{\partial y}+X_{3}%
\frac{\partial }{\partial z}+X_{4}\frac{\partial }{\partial t}+X_{5}\frac{%
\partial }{\partial s}.
\]%
Then $\left( \phi ,\xi ,\eta ,\widetilde{g}\right) $ is a spacelike almost $%
\left( \varepsilon \right) $-paracontact structure on ${\Bbb R}^{5}$. We
note that {\rm index}$\left( \widetilde{g}\right) =2$.

Now consider a hypersurface $M$ given by
\[
t=z.
\]%
Then the tangent bundle $T\!M$ of $M$ is spanned by
\[
\left \{ U_{1}=\frac{\partial }{\partial x},~U_{2}=\frac{\partial }{\partial
y},~U_{3}=\frac{\partial }{\partial z}+\frac{\partial }{\partial t},~U_{4}=%
\frac{\partial }{\partial s}\right \} ,
\]%
and $Rad~T\!M$ is spanned by $E=U_{3}$. Also the lightlike transversal
vector bundle is
\[
ltr~T\!M={\rm Span}\left \{ N=2\left( -\frac{\partial }{\partial z}+\frac{%
\partial }{\partial t}\right) \right \} .
\]%
Furthermore%
\[
\phi E=-\left( \frac{\partial }{\partial x}+\frac{\partial }{\partial y}%
\right) +(z+t)\frac{\partial }{\partial s}\in \Gamma \left( S(T\!M)\right) ,
\]

\[
\phi N=2\left( \frac{\partial }{\partial x}-\frac{\partial }{\partial y}%
+(t-z)\frac{\partial }{\partial s}\right) \in \Gamma \left( S(T\!M)\right) .
\]%
Thus $M$ is a screen semi-invariant lightlike hypersurface of ${\Bbb R}^{5}$.
\end{ex-new}

Let $M$ be a screen semi-invariant lightlike hypersurface of a $\left(
n+2\right) $-dimensional $\left( \varepsilon \right) $-almost paracontact
metric manifold $\widetilde{M}$. We set
\begin{equation}
V=\phi E\qquad {\rm and}\qquad U=\phi N.  \label{eq-SI-1}
\end{equation}%
Then, from the second equation of (\ref{eq-null-5}) and (\ref{eq-metric-1}),
we obtain
\begin{equation}
g\left( V,U\right) =1.  \label{eq-SI-2}
\end{equation}%
Therefore $\left \langle V\right \rangle \oplus \left \langle
U\right
\rangle $ is a non-degenerate vector subbundle of $S(T\!M)$ of rank
$2$. Since $\xi $ belongs to $S(T\!M)$ and
\[
g\left( V,\xi \right) =g\left( U,\xi \right) =0,
\]%
therefore there exists a non-degenerate distribution $D_{0}$ of rank $n-3$
on $M$ such that
\begin{equation}
S(T\!M)=D_{0}\perp \left \{ \left \langle V\right \rangle \oplus \left
\langle U\right \rangle \right \} \perp \left \langle \xi \right \rangle .
\label{eq-SI-3}
\end{equation}%
We note that $D_{0}$ is an invariant distribution with respect to $\phi $,
that is, $\phi D_{0}=D_{0}$. Denoting
\[
D=D_{0}\perp Rad~T\!M\perp \left \langle V\right \rangle \qquad {\rm and}%
\qquad D^{\prime }=\left \langle U\right \rangle ,
\]%
we have
\begin{equation}
T\!M=D\oplus D^{\prime }\perp \left \langle \xi \right \rangle .
\label{eq-SI-4}
\end{equation}%
Thus, every $X\in \Gamma (T\!M)$ can be expressed as
\[
X=RX+QX+\eta \left( X\right) \xi ,
\]%
where $R$ and $Q$ are the projections of $T\!M$ into $D$ and $D^{\prime }$,
respectively. Hence, we may write
\[
\varphi X=\phi RX,\qquad X\in \Gamma (T\!M).
\]%
From (\ref{eq-eta-xi}), (\ref{eq-llh-1}) and (\ref{eq-llh-2}), we obtain%
\begin{equation}
\phi ^{2}X=\varphi ^{2}X+u\left( X\right) U+u(\varphi X)N.  \label{eq-SI-5}
\end{equation}%
By comparing the tangential and transversal parts in (\ref{eq-SI-5}) we get
\begin{equation}
\varphi ^{2}=I-\eta \otimes \xi -u\otimes U,  \label{eq-SI-6}
\end{equation}%
\begin{equation}
u\circ \varphi =0,  \label{eq-SI-7}
\end{equation}%
respectively. Next, from (\ref{eq-phi-xi}) one can easily see that
\begin{equation}
\varphi \xi =0{\quad {\rm and}\quad }u\left( \xi \right) =0.  \label{eq-SI-8}
\end{equation}%
Since $\phi ^{2}N=N$, by using (\ref{eq-llh-1}) we also have%
\begin{equation}
\varphi U=0{\quad {\rm and}\quad }u\left( U\right) =1.  \label{eq-SI-9}
\end{equation}%
Furthermore, from (\ref{eq-eta-phi}) we have%
\begin{equation}
\eta \left( U\right) =0.  \label{eq-SI-10}
\end{equation}%
Finally, we get%
\[
(\eta \circ \varphi )X=\eta (\phi X-u(X)N),
\]%
which gives%
\begin{equation}
\eta \circ \varphi =0.  \label{eq-SI-11}
\end{equation}%
Thus we have the following

\begin{prop}
Let $M$ be a screen semi-invariant lightlike hypersurface of an $%
(\varepsilon )$-almost paracontact metric manifold $(\widetilde{M},\phi ,\xi
,\eta ,\overline{g},\varepsilon )$. Then $M$ possesses a para $\left(
\varphi ,\xi ,\eta ,U,u\right) $-structure, that is,
\[
\varphi ^{2}=I-\eta \otimes \xi -u\otimes U,\qquad \varphi \xi =0,\qquad
\varphi U=0,\qquad \eta \circ \varphi =0,
\]%
\[
u\circ \varphi =0,\quad \eta \left( \xi \right) =1,\quad u\left( U\right)
=1,\qquad \eta \left( U\right) =0,\qquad u\left( \xi \right) =0.
\]
\end{prop}

Next, we have the following:

\begin{th}
Let $M$ be a screen semi-invariant lightlike hypersurface of an $%
(\varepsilon )$-para Sasakian manifold $(\widetilde{M},\phi ,\xi ,\eta ,%
\overline{g},\varepsilon )$. Then we have
\begin{eqnarray}
\left( \nabla _{X}\varphi \right) Y &=&u\left( Y\right) A_{N}X+B\left(
X,Y\right) U  \nonumber \\
&&-\,g(X,Y)\xi +2\varepsilon \eta \left( X\right) \eta \left( Y\right) \xi
-\varepsilon \eta \left( Y\right) X,  \label{eq-SI-del-phi-1}
\end{eqnarray}%
\begin{eqnarray}
\left( \nabla _{X}\varphi \right) Y &=&u\left( Y\right) A_{N}X+B\left(
X,Y\right) U  \nonumber \\
&&-\, \left( g\left( \varphi X,\varphi Y\right) +u(X)\theta (\varphi
Y)+u(Y)\theta (\varphi X)\right) \xi  \nonumber \\
&&-\, \varepsilon \eta \left( Y\right) \left( \varphi ^{2}X+u\left( X\right)
U\right) ,  \label{eq-SI-del-phi-2}
\end{eqnarray}%
\begin{equation}
(\nabla _{X}u)Y=-\,B(X,\varphi Y)-u(Y)\tau (X),  \label{eq-SI-del-u}
\end{equation}%
\begin{equation}
\nabla _{X}U=-\varphi (A_{N}X)+\tau \left( X\right) U,  \label{eq-SI-del-U}
\end{equation}%
\begin{equation}
B(X,U)=-u(A_{N}X)  \label{eq-SI-B(X,U)}
\end{equation}%
for all $X,Y\in \Gamma (T\!M)$.
\end{th}

\noindent {\bf Proof.} We have
\begin{eqnarray}
(\widetilde{\nabla }_{X}\phi )Y &=&\left( \nabla _{X}\varphi \right)
Y-u\left( Y\right) A_{N}X-B\left( X,Y\right) U  \nonumber \\
&&+\left \{ \left( \nabla _{X}u\right) Y+u\left( Y\right) \tau \left(
X\right) +B\left( X,\varphi Y\right) \right \} N,  \label{eq-SI-del-phi-3}
\end{eqnarray}%
where (\ref{eq-llh-1}), (\ref{eq-null-7}), (\ref{eq-null-8}) and (\ref%
{eq-SI-1}) are used. Next, from (\ref{para2}) we have
\begin{equation}
(\widetilde{\nabla }_{X}\phi )Y=-\,g(X,Y)\xi +2\varepsilon \eta \left(
X\right) \eta \left( Y\right) \xi -\varepsilon \eta \left( Y\right) X,\qquad
X,Y\in \Gamma (T\!M).  \label{eq-SI-del-phi-4}
\end{equation}%
Using (\ref{eq-g(vphiX,vphiY)}), (\ref{eq-SI-5}) and (\ref{eq-SI-7}) in (\ref%
{para2}), for all $X,Y\in \Gamma (T\!M)$, we also have
\begin{eqnarray}
(\widetilde{\nabla }_{X}\phi )Y &=&-\, \left( g\left( \varphi X,\varphi
Y\right) +u(X)\theta (\varphi Y)+u(Y)\theta (\varphi X)\right) \xi  \nonumber
\\
&&-\varepsilon \eta \left( Y\right) \left( \varphi ^{2}X+u\left( X\right)
U\right) .  \label{eq-SI-del-phi-5}
\end{eqnarray}%
From (\ref{eq-SI-del-phi-3}) and (\ref{eq-SI-del-phi-4}) we get (\ref%
{eq-SI-del-phi-1}). Similarly, from (\ref{eq-SI-del-phi-3}) and (\ref%
{eq-SI-del-phi-5}) we get (\ref{eq-SI-del-phi-2}). Next, taking transversal
part in (\ref{eq-SI-del-phi-3}) to be zero, we get (\ref{eq-SI-del-u}).

\smallskip

Using (\ref{eq-null-7}), (\ref{eq-null-8}) and (\ref{eq-llh-1}) we get
\begin{eqnarray}
(\widetilde{\nabla }_{X}\phi )N &=&\nabla _{X}U+\varphi (A_{N}X)-\tau \left(
X\right) U  \nonumber \\
&&+(B(X,U)+u(A_{N}X))N.  \label{eq-SI-del-phi-N}
\end{eqnarray}%
Since, from (\ref{para2}) we have
\begin{equation}
(\widetilde{\nabla }_{X}\phi )N=0,  \label{eq-SI-del-tilde-phi-N}
\end{equation}%
then using (\ref{eq-SI-del-phi-N}) and \ref{eq-SI-del-tilde-phi-N} we have
\begin{eqnarray*}
0 &=&\nabla _{X}U+\varphi (A_{N}X)-\tau \left( X\right) U \\
&&+(B(X,U)+u(A_{N}X))N.
\end{eqnarray*}%
By equating the tangential and transversal parts in the previous equation we
get (\ref{eq-SI-del-U}) and (\ref{eq-SI-B(X,U)}), respectively. $%
\blacksquare $

\begin{prop}
Let $M$ be a screen semi-invariant lightlike hypersurface of an $%
(\varepsilon )$-para Sasakian metric manifold $(\widetilde{M},\phi ,\xi
,\eta ,\overline{g},\varepsilon )$. Then for any $X,Y\in \Gamma (T\!M)$, the
Lie derivative of $g$ with respect to the vector field $V$ is given by%
\begin{equation}
({\cal L}_{V}g)=X(u(Y))+Y(u(X))+u([X,Y])-2u(\nabla _{X}Y).
\label{eq-Lie-V-g}
\end{equation}
\end{prop}

\noindent {\bf Proof. }We have{\bf \ }%
\begin{equation}
({\cal L}_{Z}\widetilde{g})(X,Y)=({\cal L}_{Z}g)(X,Y),\qquad X,Y,Z\in \Gamma
(T\!M).  \label{eq-Lie-V-g-1}
\end{equation}%
Hence, using $V=\widetilde{\varphi }E$ and $\widetilde{g}((\widetilde{\nabla
}_{X}\widetilde{\varphi })E,Y)=0$ we have
\begin{equation}
({\cal L}_{V}g)(X,Y)=({\cal L}_{V}\widetilde{g})(X,Y)=\widetilde{g}(%
\widetilde{\nabla }_{X}E,\widetilde{\varphi }Y)+\widetilde{g}(\widetilde{%
\varphi }X,\widetilde{\nabla }_{Y}E).  \label{eq-Lie-V-g-2}
\end{equation}%
Next,%
\begin{eqnarray*}
\widetilde{g}(\widetilde{\nabla }_{X}E,\phi Y) &=&\widetilde{g}(\widetilde{%
\nabla }_{X}E,\varphi Y+u(Y)N) \\
&=&\widetilde{g}(\widetilde{\nabla }_{X}E,\varphi Y)+u(Y)\widetilde{g}(%
\widetilde{\nabla }_{X}E,N) \\
&=&-\widetilde{g}(E,\widetilde{\nabla }_{X}\varphi Y)-u(Y)\widetilde{g}(E,%
\widetilde{\nabla }_{X}N) \\
&=&-B(X,\varphi Y)-u(Y)\tau (X) \\
&=&(\nabla _{X}u)Y,
\end{eqnarray*}%
which gives%
\begin{equation}
\widetilde{g}(\widetilde{\nabla }_{X}E,\phi Y)=X(u(Y))-u(\nabla _{X}Y).
\label{eq-Lie-V-g-3}
\end{equation}%
Using (\ref{eq-Lie-V-g-3}) in (\ref{eq-Lie-V-g-2}) we complete the proof. $%
\blacksquare $

\subsection{Integrability of $D\perp \left \langle \protect\xi %
\right
\rangle $}

We note that $X\in \Gamma (D\perp \left \langle \xi \right \rangle )$ if and
only if $u\left( X\right) =0$. Now from (\ref{eq-SI-del-u}), we have for all
$X,Y\in \Gamma (TM)$
\[
u(\nabla _{X}Y)=\nabla _{X}u(Y)+B(X,\varphi Y)+u(Y)\tau (X),
\]%
from which we get
\begin{eqnarray*}
u\left[ X,Y\right] &=&B(X,\varphi Y)-B(\varphi X,Y) \\
&&+\nabla _{X}u(Y)-\nabla _{Y}u(X)+u(Y)\tau (X)-u(X)\tau (Y).
\end{eqnarray*}%
Now, let $X,Y\in \Gamma (D\perp \left \langle \xi \right \rangle )$. Then $%
u\left( X\right) =0=u\left( Y\right) $, and from the previous eqiation we
get
\[
u\left[ X,Y\right] =B(X,\varphi Y)-B(\varphi X,Y)
\]%
for all $X,Y\in \Gamma (D\perp \left \langle \xi \right \rangle )$. Thus we
get a necessary and sufficient condition for the integrability of the
distribution $D\perp \left \langle \xi \right \rangle $ in the following:

\begin{th}
Let $M$ be a screen semi-invariant lightlike hypersurface of an $%
(\varepsilon )$-para Sasakian manifold $(\widetilde{M},\phi ,\xi ,\eta ,%
\overline{g},\varepsilon )$. Then\ the distribution\ $D\perp \left \langle
\xi \right \rangle $ is integrable if and only if
\[
B(X,\varphi Y)=B(\varphi X,Y),\qquad X,Y\in \Gamma (D).
\]
\end{th}

\subsection{Integrability of $D^{\prime }\perp \left \langle \protect\xi %
\right \rangle $}

Here we find a necessary and sufficient condition for the distribution $%
D^{\prime }\perp \left \langle \xi \right \rangle $ to be integrable.

\begin{th}
Let $M$ be a screen semi-invariant lightlike hypersurface of an $%
(\varepsilon )$-para Sasakian manifold $(\widetilde{M},\phi ,\xi ,\eta ,%
\overline{g},\varepsilon )$. Then\ the distribution\ $D^{\prime }\perp
\left
\langle \xi \right \rangle $ is integrable if and only if
\begin{equation}
A_{N}\xi +\varepsilon U=0.  \label{eq-SI-D'-xi-int}
\end{equation}
\end{th}

\noindent {\bf Proof.} Note that $X\in \Gamma (D^{\prime }\perp
\left
\langle \xi \right \rangle )$ if and only if $\varphi X=0$. Now for
all $X,Y\in \Gamma (TM)$, in view of (\ref{eq-SI-del-phi-1}), we have
\begin{eqnarray*}
\varphi \left( \nabla _{X}Y\right) &=&\nabla _{X}\varphi Y-u\left( Y\right)
A_{N}X-B\left( X,Y\right) U \\
&&+\,g(X,Y)\xi -2\varepsilon \eta \left( X\right) \eta \left( Y\right) \xi
+\varepsilon \eta \left( Y\right) X.
\end{eqnarray*}%
From the above equation we get
\begin{eqnarray*}
\varphi \left[ X,Y\right] &=&\nabla _{X}\varphi Y-\nabla _{Y}\varphi
X+u\left( X\right) A_{N}Y-u\left( Y\right) A_{N}X \\
&&+\varepsilon \eta \left( Y\right) X-\varepsilon \eta \left( X\right) Y.
\end{eqnarray*}%
In particular, for $X,Y\in \Gamma (D^{\prime }\perp \left \langle \xi
\right
\rangle )$ we get
\begin{equation}
\varphi \left[ X,Y\right] =u\left( X\right) A_{N}Y-u\left( Y\right)
A_{N}X+\varepsilon \eta \left( Y\right) X-\varepsilon \eta \left( X\right) Y.
\label{eq-SI-int-1}
\end{equation}%
But $D^{\prime }$ and $\left \langle \xi \right \rangle $ are integrable,
hence $D^{\prime }\perp \left \langle \xi \right \rangle $ is integrable if
and only if
\[
\varphi \left[ U,\xi \right] =0,
\]%
which, in view of (\ref{eq-SI-int-1}), is equivalent to (\ref%
{eq-SI-D'-xi-int}). $\blacksquare $


\begin{thebibliography}{99}
\bibitem{Beem-Ehrl-81} J.K. Beem and P.E. Ehrlich, {\em Global Lorentzian
geometry}, Marcel Dekker, New York, 1981.

\bibitem{Bej-Dug-93} A. Bejancu and K.L. Duggal, {\em Real hypersurfaces of
indefinite Kaehler manifolds}, Internat. J. Math. Math. Sci. {\bf 16}
(1993), no. 3, 545-556.

\bibitem{Blair-02-book} D.E. Blair, {\em Riemannian geometry of contact and
symplectic manifolds,} Progress in Mathematics, 203. Birkhauser Boston,
Inc., Boston, MA, 2002.

\bibitem{Dug-90-IJMMS} K.L. Duggal, {\em Space time manifolds and contact
structures}, Internat. J. Math. Math. Sci. 13 (1990), no. 3, 545--553.

\bibitem{Dug-Bej-96} K.L. Duggal and A. Bejancu, {\em Lightlike submanifolds
of semi-Riemannian manifolds and its applications}, Kluwer, Dordrecht 1996.

\bibitem{Dug-Jin-07} K.L. Duggal and D.H. Jin, {\em Null Curves and
Hypersurfaces of Semi-Riemannian Manifolds}, World Scientific Publishing Co.
Pvt. Ltd., 2007.

\bibitem{Dug-Sah-07} K.L. Duggal and B. Sahin, {\em Lightlike submanifolds
of indefinite Sasakian manifolds}, Int. J. Math. Math. Sci. 2007, Art. ID
57585, 21 pp.

\bibitem{Dug-Sah-10-book} K.L. Duggal and B. Sahin, {\em Differential
geometry of lightlike submanifolds}, Birkh\"{a}user, 2010.

\bibitem{Kang-JKP-03} T.H. Kang, S.D. Jung, B.H. Kim, H.K. Pak and J.S. Pak,
{\em Lightlike hypersurfaces of indefinite Sasakian manifolds}, Indian J.
Pure Appl. Math. {\bf 34} (2003), no. 9, 1369-1380.

\bibitem{Mat-89} K. Matsumoto, {\em On Lorentzian paracontact manifolds},
Bull. Yamagata Univ. Natur. Sci. {\bf 12} (1989), no. 2, 151-156.

\bibitem{ONeill-83} B. O'Neill, {\em Semi-Riemannian geometry with
applications to relativity}, Academic Press, 1983.

\bibitem{Sasaki-60-Tohoku} S. Sasaki, {\em On differentiable manifolds with
certain structures which are closely related to almost contact structure I},
T\^{o}hoku Math. J. {\bf 12} (1960), 459-476.

\bibitem{Sato-76} I. Sat\={o}, {\em On a structure similar to the almost
contact structure}, Tensor (N.S.) {\bf 30} (1976), no. 3, 219-224.

\bibitem{Takahashi-69-Tohoku-1} T. Takahashi, {\em Sasakian manifold with
pseudo-Riemannian metric}, T\^{o}hoku Math. J. {\bf 21} (1969), 644-653.

\bibitem{Tri-KYK-10} M.M. Tripathi, E. K\i l\i \c{c}, S. Y\"{u}ksel Perkta%
\c{s} and S. Kele\c{s}, {\em Indefinite almost paracontact metric manifolds}%
, Int. J. Math. Math. Sci. 2010 (2010), Art. Id. 846195, pp. 19.

\bibitem{Yuk-KTK-12} S. Y\"{u}ksel Perkta\c{s}, E. K\i l\i \c{c}, M.M.
Tripathi and S. Kele\c{s}, {\em On} $(\varepsilon )$-{\em para Sasakian} $3$%
{\em -manifolds}, Int. J. Pure Appl. Math. {\bf 77} (2012), no.4, 485-499.
\end{thebibliography}
\end{document}